\documentclass[a4paper,11pt]{amsart}
\usepackage{amsfonts,amssymb}
\author{O. V. Chuvashova}
\address{
Department of Higher Algebra\\
Faculty of Mechanics and Mathematics\\
Moscow State University\\
119992 Moscow, Russia} \email{chuvashova@yandex.ru}

\title[The separation properties]{The separation properties for closures of toric orbits}

\date{September 23, 2004}



\newtheorem{theorem}{Theorem}
\theoremstyle{conjecture}

\theoremstyle{definition}
\newtheorem{definition}{Definition}
\newtheorem{example}{Example}
\theoremstyle{proposition}
\newtheorem{proposition}{Proposition}
\theoremstyle{lemma}
\newtheorem{lemma}{Lemma}
\theoremstyle{corollary}
\newtheorem{corollary}{Corollary}
\theoremstyle{remark}
\newtheorem{remark}{Remark}
\newcommand{\theo}{\par \noindent
                   {\bf Theorem}
                      }

\newcommand{\g}{\mathfrak{g}}
\newcommand{\ssl}{\mathfrak{sl}}

\newcommand{\e}{\epsilon}
\newcommand{\al}{\alpha}
\newcommand{\be}{\beta}
\newcommand{\ga}{\gamma}
\newcommand{\PProj}{{\mathbb{ P}}(V)}
\newcommand{\Pd}{{\mathbb{P}}(V^*)}
\newcommand{\N}{\mathbb{N}}

\newcommand{\Q}{\mathbb{Q}}
\newcommand{\Z}{\mathbb{Z}}
\newcommand{\la}{\langle}
\newcommand{\ra}{\rangle}

\newcommand{\codim}{\rm{codim}}

\begin{document}

\begin{abstract}

A subset $X$ of a vector space $V$ is said to have the "Separation
Property" if it separates linear forms in the following sense:
given a pair $(\al, \be)$ of linearly independent forms on $V$
there is a point $x$ on $X$ such that $\al(x)=0$ and $\be(x)\ne
0$. A more geometric way to express this is the following: every
homogeneous hyperplane $H\subseteq V$ is linearly spanned by its
intersection with $X$. The separation property was first asked for
conjugacy classes in simple Lie algebras.

We give an answer for orbit closures in representation spaces of
an algebraic torus. We consider also the strong and the weak
separation properties. It turns out that toric orbits well
illustrate these concepts.

\end{abstract}

\maketitle

\section{Introduction}

1. Let $V$ be a vector space over a field $k$ and $G$ be a
connected reductive algebraic group. Consider a linear action
$G:V$. Properties of the closure of an orbit $\overline{Gv}\subset
V$ can be sorted into four groups:

(1) "combinatorial" (the number of orbits in $\overline{Gv}$, the
graph of the orbit adherence relation, $\ldots$),

(2) geometrical (the smoothness, the normality, the Cohen-Macaulay
property, types of singularities, $\ldots$),

(3) topological (the contractibility, the simple connectedness,
the computation of homologies and cohomologies, higher homotopy
groups, $\ldots$),

(4) properties of the embedding $\overline{Gv}\subset V$ (the
dimension of the linear span, hyperplane sections, a description
of the ideal defining the variety,~$\ldots$).

In this paper we consider the separation properties. In our
opinion they belong to the most natural properties of the fourth
group.

\begin{definition} A subset $X$ of a vector space $V$ has the {\it
separation property} (briefly (SP)) if for any pair of linearly
independent linear functions $\alpha, \beta \in V^*$ there  exists
a point $x\in X$ such that $\al (x)=0$ and $\be (x)\ne 0$.
\end{definition}

In \ other words, \ the separation \ property \ for $X\subset V$ \
means \ that \ $H\cap X\not\subseteq H'$ for any pair $H\ne H'$ of
homogeneous hyperplanes in $V$. Or, equivalently, for any
homogeneous hyperplane $H$ the intersection $H\cap X$ linearly
spans $H$.

For projective spaces we have a similar definition: a subset
$Y\subset \PProj$ has the separation property if for any pair
$D\ne D'$ of hyperplanes in $\PProj$ we have $D\cap Y\not\subseteq
D'$.

\begin{remark} \ Let $\al\in V^*$ and $H_{\al}$ be the corresponding
hyperplane. If either $X\cap H_{\al}$ is empty or $X\cap H_{\al}$
is zero or $X\subseteq H_{\al}$, then it is easy to see that $X$
has not (SP).\end{remark}

For the first time the separation property was asked by
J.~C.~Jantzen in connection with the paper of A.~Premet
\cite{Prem}:

\medskip

\noindent {\bf Question.} \ Let $k$ be an algebraically closed
field, $G$ be a simple algebraic group, and $\g$ be its tangent
algebra. Is it true that the minimal nilpotent orbit of the
adjoint action in $\g$ has the separation property?

\medskip

The answer was obtained in the work \cite{KrW}. It is affirmative
for all simple groups except for $Sp_{2n}$.

\begin{example} Consider the Lie algebra $\ssl_2$. Here the minimal
nilpotent orbit consists of the following matrices:

$$
M=\left\{ \left(
\begin{array}{cc}
a & b \\
c & -a
\end{array}
\right) : a^2+bc=0 \right\}.
$$
It is easy to see that $M$ has not the separation property. In
fact, $b=0$ implies $a=0$ for a matrix from $M$. \end{example}

The notions of the "strong"  and the "weak"  separation properties
also were introduced in \cite{KrW}.

\begin{definition} A closed affine cone $X\subseteq V$ of dimension $\geq 2$
has the {\it strong separation property} (briefly (SSP)) if for
any linear subspace $W\subseteq V$ of codimension 2 we have
$\codim_X W\cap X = 2$.\end{definition}

There is  a  similar  definition  for  a  closed  projective
subvariety. A closed subvariety $Y\subseteq \PProj$ has the strong
separation property if for any  linear  subspace  $L\subseteq
\PProj$ of  codimension  2  we  have  $\codim_Y L\cap Y=2$.

\begin{remark} The strong separation property for a closed projective
subva\-riety $Y\subseteq \PProj$ (closed affine cone $X\subseteq
V$) implies the separation property. \end{remark}

\begin{proof} Consider the projective case. Suppose that (SP) does
not hold. This means that there exist two different hyperplanes
$H_{\al}$ and $H_{\be}$ in $\PProj$ such that $H_{\al}\cap
Y\subseteq H_{\be}$. Then $(H_{\al}\cap H_{\be})\cap Y =
H_{\al}\cap Y$ has the codimension $\leq 1$ in $Y$ (see
\cite[Th.~1.6.4]{Sh}) and (SSP) does not hold.

In the affine case it is sufficient to notice that (SP) (resp.
(SSP)) for a closed affine cone is equivalent to (SP) (resp.
(SSP)) for its projectivization. \end{proof}

The next example shows that (SP) does not imply (SSP).

\begin{example} Consider the subvariety of degenerate matrices in the space
of all $2\times 2$ matrices:
$$
N=\left\{ \left(
\begin{array}{cc}
a & b \\
c & d
\end{array}
\right) : ad-bc=0 \right\}.
$$
It has not (SP) as it contains the subspace $a=b=0$. By direct
calculation it is easy to check that $N$ has (SP). (It also
follows from Theorem 3 of this work. Consider the linear action
$T=(k^*)^3 : k^4$, $(t_1,t_2,t_3)(x_1,x_2,x_3,x_4)=(t_1x_1,t_3
x_2,t_1t_2x_3,t_2t_3x_4)$. Note that $N$ is isomorphic to
$\overline{T(1,1,1,1)}$.)
\end{example}

\begin{definition} A subset $X$ of a vector space $V$ has the {\it weak
separation property} (briefly (WSP)) if for any pair of
homogeneous hyperplanes $H\ne H'$ we have $H\cap X\ne H'\cap X$
(in the set-theoretical sense).\end{definition}

The definition of the weak separation property for a subset of a
projective space is similar.

It is obvious that the separation property implies the weak
separation property.

\begin{example} It is easy to check that the subvariety $M$ of nilpotent
matrices in $\ssl_2$ has (WSP) and has not (SP). (This also
follows from Theorems 3 and 4 of this work. Consider the linear
action $T=(k^*)^2 :k^3$, $(t_1,t_2)(x_1,x_2,x_3)=(t_1t_2x_1,
t_1^2x_2, t_2^2x_3)$. Note that $M$ is isomorphic to
$\overline{T(1,1,1)}$.) \end{example}

The following theorems are proved in \cite{KrW}.

\medskip
\theo{\bf \cite[Th. 1]{KrW}.} {\it Let $G$ be a connected
semisimple group and $V(\lambda)$ be its Weyl module corresponding
to a highest weight $\lambda$ with numerical marks $n_i$. Denote
by $O_{min}\subset V(\lambda)$ the orbit of a highest vector.
Then:

(1) $\overline{O_{min}}$ satisfies (SSP) $\iff$ $\lambda$ is a
fundamental weight;

(2) $O_{min}$ satisfies  (SP) $\iff$ $n_i\leq 1$ for any $i$;

(3) $O_{min}$ satisfies (WSP) $\iff$ $n_i\leq 2$ for any $i$.}

\medskip
\theo{\bf \cite[Th. 2]{KrW}.} {\it Let $G:V$ be an irreducible
representation of a connected semisimple group $G$. If $O_{min}$
satisfies (SP), then any $G$-orbit in $V$ satisfies (SP).}

\medskip
\theo{\bf \cite[Th. 3]{KrW}.} {\it Suppose that ${\rm char}\, k =
0$ and $G:V$ is an irreducible representation of a semisimple
group $G$. Then a generic $G$-orbit in $V$ satisfies (SP).}
\medskip

2. The aim of this work is to investigate the separation
properties for closures of toric orbits in vector and projective
spaces over an algebraically closed field. This is a primary
generalization of Theorem \cite[Th.~1]{KrW} to the case of
reducible representations of reductive groups.

Let $T$ be an algebraic torus, $\Lambda$ be the lattice of
characters of $T$, and $V$ be a vector space over an algebraically
closed field $k$. Consider a linear action $T:V$, where
$$t(x_1,\ldots,x_n)=(\chi_1(t)x_1,\ldots,\chi_n(t)x_n).$$ Let
$\Sigma$ be the semigroup in $\Lambda$ generated by the characters
$\chi_1,\ldots,\chi_n$ and $K$ be the cone in
$\Lambda\otimes_{\Z}\Q$ generated by $\Sigma$.
\medskip

\theo {\bf 3.} {\it The orbit closure
$X=\overline{T(1,\ldots,1)}\subset V$ satisfies (SP) if and only
if the following conditions hold:

\noindent (1) the cone $K$ is strictly convex;

\noindent (2) $\Q_+\chi_i$ is an edge of $K$ for any $i$;

\noindent (3) $\Q_+\chi_i\ne\Q_+\chi_j$ for $i\ne j$.}
\medskip

For a projective action of a torus we have:

\medskip

\theo {\bf 3'.} {\it The orbit closure
$X=\overline{T(1:\ldots:1)}\subset \PProj$ satisfies (SP) if and
only if the following conditions hold:

\noindent (1) the point $\chi_i$ is a vertex of the convex hull
${conv}\{\chi_1,\ldots,\chi_n\}$  for any $i$;

\noindent (2) $\chi_i\ne \chi_j$ for $i\ne j$.}
\medskip

The proofs are based on the fact that if the separation property
fails on some pair of hyperplanes, then there exist such
$T$-invariant hyperplanes. To prove this statement we introduce
the notion of the characteristic variety of a subset $X\subset V$
(or $X\subset \PProj$):
$$Ch(X)=\{(\la\al\ra,\la\be\ra)\in \Pd\times\Pd \ | \ \al(x)=0
\Longrightarrow \be(x)=0 \quad\forall x\in X\}.$$

After that we prove:
\medskip
\theo {\bf 2.} {\it Suppose that an affine subvariety $X\subset V$
is irreducible, is not contained in a homogeneous hyperplane,
meets any homogeneous hyperplane and $\dim X > 1$. Then $Ch(X)$ is
closed in ${\Pd}\times {\Pd}$.}
\medskip

Finally we apply the fact that an algebraic torus acting on a
projective variety has a fixed point. More precisely:

\medskip
\noindent {\bf Proposition 3.} {\it If $X$ does not satisfy (SP)
and $Ch(X)$ is closed, then there exists a pair
$(\la\al\ra,\la\be\ra)\in Ch(X)$ such that $\al, \be$ are
$T$-semiinvariant and linearly independent.}
\medskip

Proposition 3 allows to simplify the proof of the criterion of the
separation property for $SL_2$-orbits of binary forms obtained in
the thesis of K.~Baur.

\medskip
\theo{\bf \cite[Th. 3.4]{Baur}.} {\it Let $f\in k[x,y]_n$. Then
the orbit $O_f=SL_2 f$ satisfies (SP) if and only if  $f$ has a
linear factor of multiplicity one.}

\medskip

For the weak separation property we obtain the following theorems.

\medskip
\theo {\bf 4.} {\it The orbit closure
$X=\overline{T(1,\ldots,1)}\subset V$ satisfies (WSP) if and only
if the following conditions hold:

\noindent (1) the cone $K$ is strictly convex;

\noindent (2) there are no more then one characters $\chi_i$ in
the interior of any face of $K$ (in particular,
$\Q_+\chi_i\ne\Q_+\chi_j$ for $i\ne j$).}

\medskip
\theo {\bf 4'.} {\it The orbit closure
$X=\overline{T(1:\ldots:1)}\subset \PProj$ satisfies (WSP) if and
only if  there are no more then one characters $\chi_i$ in the
interior of any face of the convex hull
$conv\{\chi_1,\ldots,\chi_n\}$ (in particular, $\chi_i\ne\chi_j$
for $i\ne j$).}

\medskip

Finally we consider the strong separation property.

\medskip
\theo {\bf 5.} {\it Suppose that the orbit closure
$X=\overline{T(1,\ldots,1)}\subset V$ is a cone; then $X$
satisfies (SSP) if and only if $X=V$ (i.e., the weights
$\chi_1,\ldots,\chi_n$ are linearly independent).}

\medskip
\theo {\bf 5'.} {\it The orbit closure
$X=\overline{T(1:\ldots:1)}\subset \PProj$  satisfies (SSP) if and
only if $X=\PProj$ (i.e., the weights $\chi_1,\ldots,\chi_n$ are
affinely independent).}
\medskip

I am grateful to my adviser Ivan V. Arzhantsev for the subject of
this work, the permanent support, and numerous remarks and ideas.
I also thank Dmitri A. Timashev for important comments and
detection of a gap in the preliminary version of the text.

\section{Hypersurfaces}

\begin{proposition} {\it Let $S\subset V$ be a hypersurface. Then

(1) $S$ does not satisfy (SP) $\iff$ $S$ can be defined by an
equation  $x_2^n + x_1F=0$  in some coordinate system for some
$n\in \N$ and $F\in k[V]$;

(2) $S$ does not satisfy (WSP) $\iff$ $S$ can be defined by an
equation $x_1^n + x_2^m + x_1x_2F=0$ in some coordinate system for
some $n, m\in \N$ and $F\in k[V]$;

(3) a homogeneous hypersurface $S$ does not satisfy (SSP) $\iff$
 $S$ can be defined by an equation
$x_1F_1 + x_2F_2=0$ in some coordinate system for some homogeneous
polynomials $F_1, F_2\in k[V]$ of the same degree.}
\end{proposition}

\begin{proof} Let $S$ be defined by an equation $P=0$. We may assume
 that $P$ has no multiple factors. Since the field $k$ is
algebraically closed, it follows that $P$ is defined uniquely up
to a constant.

(1)\ \ If\ \ $P$\ \ has\ \ the\ \ form\ \ $x_2^n + x_1F$\ \ in\ \
some\ \  coordinate\ \ system,\ \ then\ \ $S\cap H_{x_1} \subseteq
H_{x_2}$ and $S$ has not (SP).

Conversely, if $S$ has not (SP), then we can choose a coordinate
system such that $S\cap H_{x_1} \subseteq H_{x_2}$.  This  means
that  $P(0,x_2,\ldots)=0$  implies  $x_2=0$. By Hilbert's
nullstellentsatz, we have $x_2^l=P(0,x_2,\ldots)f(x_2,\ldots)$.
This implies $P(0,x_2,\ldots)=cx_2^n$ for some $n\in \N, n\leq l$,
and $c\in k^*$. Then $P = cx_2^n + x_1F$. (We may assume that
$c=1$.)

(2) If $P$ has the form $x_1^n + x_2^m + x_1x_2F$ in some
coordinate system, then $S\cap H_{x_2} = S\cap H_{x_1}$ and $S$
has not (WSP).

Conversely, if $S$ has not (WSP), then we can choose a coordinate
system such that $S\cap H_{x_2} = S\cap H_{x_1}$. Statement (1)
implies that in this coordinate system $P = x_1^n + x_2F_1 = x_2^m
+ x_1F_2$ (up  to the proportionality of the basis vectors). This
implies $P=x_1^n + x_2^m + x_1x_2F$.

(3) If $P$ has the form $x_1F_1 + x_2F_2$ in some coordinate
system, then the subspace in $V$ defined by the equations $x_1=0$,
$x_2=0$ has the codimension one in $S$ and $S$ has not (SSP).

Conversely, if $S$ has not (SSP), then we can choose a coordinate
system such that the subspace in $V$ defined by the equations
$x_1=0$, $x_2=0$ has the codimension $\leq 1$ in $S$ (this implies
that this subvariety is contained in $S$). The ideal generated by
the polynomials $x_1$, $x_2$ is radical, so, by Hilbert's
nullstellentsatz, $P=x_1F_1+x_2F_2$. \end{proof}

Denote by $I(X)$ the ideal of a closed affine subvariety
$X\subseteq V$.

\begin{proposition} {\it Let $X\subseteq V$ be a closed affine subvariety. Then

(1) $X$ does not satisfy (SP) if and only  if $X$ is contained in
a hyperpsurface $S$ such that $S$ does not satisfy (SP);

(2)  $X$ does not satisfy (WSP) if and only  if $X$ is contained
in a hypersurface $S$ such that $S$ does not satisfy
(WSP).}\end{proposition}

\begin{proof} Suppose that $X$ has (SP) (resp.  (WSP)). It
is evident that any subset $Z$ in $V$ such that $X\subseteq Z$ has
(SP) (resp. (WSP)).

Now we shall prove the converse implication. Suppose that $X$ has
not (SP). \  Then\  there\  exist\ \ linearly\ \ independent\ \
$\al,\  \be \in V^*$ \  such \  that\  $H_{\al}\cap X\subseteq
H_{\be}$. By \ Hilbert's nullstellentsatz,\ \ this\ holds\ if\
and\ only\ if\ $\be^n = f + g\al$ for some $n\in \N$, $g\in k[V]$,
$f\in I(X)$. Let $S$ be the zero set of $f$. Then $X \subseteq S$
and, by Proposition 1, $S$ has not  (SP).

Suppose that $X$ has not  (WSP). Then there exist linearly
independent $\al, \be \in V^*$ such that $H_{\al}\cap X =
H_{\be}\cap X$. Let $F$ be the hypersurface in $V$ defined by the
equation $\al(x)\be(x)=0$. The equality $H_{\al}\cap X =
H_{\be}\cap X$ is equivalent to two equalities: $H_{\al}\cap X =
F\cap X$ and $H_{\be}\cap X = F\cap X$. They hold if and only if
$rad(I(X), \al) = rad(I(X),\al\be)$ and $rad(I(X), \be) =
rad(I(X),\al\be)$ (here $(I(X),f)$ is the ideal generated by
$I(X)$ and $f$). We have $$\al^k = f_1 + g_1\al\be ,$$
$$\be^l = f_2 + g_2\al\be $$ for some $k, l\in \N$, $f_1, f_2\in I(X)$, $g_1, g_2\in
k[V]$. Summing these equations, we get $f_1+f_2 =
\al^k+\be^l-\al\be(g_1+g_2)$. Let $S$ be the zero set of the
polynomial $f_1 + f_2$. Then $S$ contains $X$ and, by the previous
proposition, has not (WSP). \end{proof}

\begin{remark} Let us remark that the analogous statement
does not hold for (SSP). In fact, let $X$ be a closed affine cone
of dimension $\geq 2$ and $f\in I(X)$ be a homogeneous polynomial.
Then the homogeneous hypersurface defined by the equation
$x_1f(x)=0$ contains $X$ and has not (SSP). Moreover, we shall
give an example of a closed affine cone having (SSP) and
containing in an irreducible hypersurface such that this
hypersurface has not (SSP).\end{remark}

\begin{example} Let $X\subset V$ be a closed irreducible affine cone of
codimension $\geq 2$ having (SSP). There exist homogeneous coprime
$f,p\in I(X)$. Consider the subspace $W=V\oplus k^2$ ($x,y$ are
coordinates in the last two items). Denote by $S$ the hypersurface
in $W$ defined by the equation $fx^m+py=0$ (where $\deg p = \deg f
+ m - 1$) and by $Y$ the closed affine cone $X\oplus k^2\subset
W$. Then $S$ has  not (SSP) and is irreducible. Note that
$Y\subset S$.

We shall prove that $Y$ has (SSP). It is sufficient to prove that
$X\subset V$ has (SSP) implies $X'=X\oplus k\subset V'=V\oplus k$
has (SSP). Assume the converse. Let (SSP) fail for $X'$ on the
subspace $U'\subset V'$ defined by equations $\al'=\al+az=0,
\be'=\be+bz=0$, where $z$ is a coordinate in the second item of
$V\oplus k$.

The first case. Let $\al$ and $\be$ be not proportional. Here
$$\dim (X\cap H_{\al}\cap H_{\be})=\dim (X'\cap H_{\al'}\cap
H_{\be'}\cap H_z)\geq \dim X'- 2=\dim X-1$$ and we obtain a
contradiction with the fact that $X$ has  (SSP).

The second case. If $\al=c\be$, then $U'$ can be defined by the
equations $\al=0$, $z=0$. As $(X'\cap H_z)\cap H_{\al}=X\cap
H_{\al}$ has the codimension 2 in $X'$, we have a
contradiction.\end{example}

\noindent {\bf Proposition 2'.} {\it Let $Y\subseteq \PProj$ be a
closed projective subvariety. Then

(1) $Y$ does not satisfy (SP) if and only if $Y$ is contained in a
hypersurface $R$ such that it does not satisfy (SP);

(2) $Y$ does not satisfy (WSP) if and only if $Y$ is contained in
a hypersurface $R$ such that it does not satisfy (WSP).}

\begin{proof} Consider the vector space $V$ and the closed affine cone
$X\subseteq V$ over $Y$. The cone $X$ has  (SP) (resp. (WSP)) if
and only if $Y$ has  (SP) (resp. (WSP)).

(1) If $Y\subseteq R$ and $R$ has not (SP), then $X\subseteq S$,
where $S$ is the affine cone over $R$ in $V$ and $S$ has not (SP).
This implies that $X$ has not (SP).

Conversely, if $X$ has not  (SP), then, by Proposition 1, there
exists a hypersurface $S\subset V$ such that $X\subseteq S$ and
$S$ has not (SP). Let $S$ be defined by the equation $f=0$, where
$f\in I(X)$. By Proposition 1, $f$ has the form $x_2^n + x_1F$ in
some coordinate system. As the ideal $I$ is homogeneous, the
homogeneous component of degree $n$ of this equation belongs to
$I(X)$. This component has the form $x_2^n + x_1F'$. Hence the
corresponding projective hypersurface contains $Y$ and has not
(SP).

(2) If $Y\subseteq R$ and $R$ has not (WSP), then $X\subseteq S$,
where $S$ is the affine cone over $R$ in $V$ and $S$ has not
(WSP). This implies that $X$ has not  (WSP).

If $X$ has not  (WSP), then, by Proposition 1, there exists a
hypersurface $S\subset V$ such that $X\subseteq S$ and $S$ has not
(WSP). Let $S$ be defined by $f=0$, where $f\in I(X)$. By
Proposition 1, $f$ has the form $x_1^n + x_2^m + x_1x_2F$ in some
coordinate system. As the ideal $I(X)$ is homogeneous, the
homogeneous components of degrees $n$ and $m$ of this equation
belong to $I(X)$. If $n=m$, then this component has the form
$x_1^n + x_2^n + x_1x_2F'$ and corresponding projective
hypersurface contains $Y$ and has not (WSP). If $n\ne m$, then
these components have the forms $x_1^n + x_1x_2F_1$ and $x_2^m +
x_1x_2F_2$. Consider the polynomial $(x_1^n + x_1x_2F_1)^m +
(x_2^m + x_1x_2F_2)^n$. This homogeneous polynomial of degree $mn$
belongs to $I(X)$ and the projective hypersurface defined by it
has not (WSP). \end{proof}

It was proved in \cite{KrW} that the strong separation property is
a property of open type in the following sense. Recall that a
family of $d$-dimensional subvarieties in $\PProj$ is a closed
subvariety $F\subset B\times \PProj$, where $B$ is an algebraic
variety such that the projection $pr_B$ induces the surjective
morphism $p:F\to B $ and any fiber of this morphism has the
dimension $d$.

\medskip
\noindent {\bf Proposition \cite[Prop. 6]{KrW}.} {\it Let $p:F\to
B$ be a family of $d$-dimensional closed subvarieties in $\PProj$.
Then the subset $\{b\in B \ | \ p^{-1}(b)$ satisfies (SSP)$\}$ is
open in $B$.}
\medskip

The following example shows that the weak separation property and
the separation property are not properties of open type.

\begin{example} Consider the closed family of hypersurfaces $F\subset
\PProj\times k$ defined by the equation
$b(x_1^m+x_2^m)+x_1x_2P=0$, where $b\in k$ and the hypersurface
$R\subset \PProj$ defined by the equation $P=0$ has  (SP) ($\deg
P=m-2$). Then, by Proposition 1, the set $p^{-1}(b)$ has not (WSP)
if $b\ne 0$. The set $p^{-1}(0)$ has  (SP), as it contains the
subset $R$ having  (SP).\end{example}

 The author does not know an example of such family of subvarieties
with an irreducible fiber over any point  $b\in B$.

\section{Characteristic varieties}

Let $X$ be a subset in a vector space $V$.

\begin{definition} The {\it characteristic variety} $Ch(X)$ of a subset $X$ is
the subset in $\Pd\times\Pd$ consisting of the following pairs:
$$Ch(X)=\{(\la\al\ra,\la\be\ra)\in \Pd\times\Pd \ | \ \al(x)=0
\Longrightarrow \be(x)=0 \quad\forall x\in X\}.$$
\end{definition}

\begin{definition} The {\it weak characteristic variety} $Chw(X)$ of a subset
$X$ is the subset in $\Pd \times \Pd$ consisting of the following
pairs:
$$Chw(X)=\{(<\al>,<\be>)\in \Pd\times\Pd \ | \
\al(x)=0 \iff $$ $$\be(x)=0 \quad\forall x\in X\}.$$
\end{definition}

\begin{remark}

\noindent (1) If $\phi:\Pd\times\Pd\to \Pd\times\Pd$ is the
symmetry with \ respect \ to \ the \ diagonal,\ i.e.,\
$\phi((\la\al\ra,\la\be\ra))=(\la\be\ra,\la\al\ra)$, then
$Chw(X)=Ch(X)\cap\phi(Ch(X))$;

\noindent (2) the diagonal $D=(\la\al\ra,\la\al\ra)\subseteq
Chw(X)\subseteq Ch(X)$;

\noindent (3) $X$ satisfies (SP) $\iff Ch(X)=D$;

\ $X$ satisfies (WSP) $\iff Chw(X)=D$;
\end{remark}

We shall need the following theorem.

\begin{theorem} \cite[Th. 4.5]{Hu} {\it Let $\phi:X\to Y$ be a dominant
morphism of irreducible varieties and $r=\dim X-\dim Y$. Suppose
that for any closed irreducible subvariety $W\subseteq Y$ any
irreducible component of $\phi^{-1}(W)$ has the dimension $\dim
W+r$. Then $\phi$ is an open morphism.} \end{theorem}

\begin{theorem} {\it Suppose that an algebraic subvariety $X\subset V$ is
irreducible, is not contained in a homogeneous hyperplane, meets
any homogeneous hyperplane and $\dim X > 1$. Then $Ch(X)$ and
$Chw(X)$ are closed.} \end{theorem}

\begin{remark} If the conditions of the theorem hold for any irreducible
component $X_i$ of a subvariety $X\subset V$, then $Ch(X)=\cap
Ch(X_i)$ is closed. \end{remark}

\begin{proof} Consider the closed subvariety
$$M=\{(\la\al\ra,\la\be\ra,x)\in\Pd\times\Pd\times X \ | \
\al(x)=0\}.$$ Let $\phi:M\to\Pd\times\Pd$ be the projection along
the third component. Consider the open subset
$L=\{(\la\al\ra,\la\be\ra,x)\in M:\be(x)\ne 0\}$ in $M$. Note that
$Ch(X)=\Pd\times\Pd\setminus \phi(L)$. Thus we need to prove that
$\phi(L)$ is open. We shall show that $\phi$ is an open morphism
applying Theorem 1.

(1) The variety $M$ is irreducible. We prove this in Lemma 1
below.

(2) The morphism $\phi$ is surjective. Indeed, $\phi(M) =
\{(\la\al\ra,\la\be\ra) \ |$ $ \exists\ x\in X : \al(x)=0\}$. But
$X$ meets any homogeneous hyperplane.

(3) Let $W\subseteq \Pd\times\Pd$ be closed and irreducible. We
have $\phi^{-1}(W)=(W\times X)\cap R$, where $R\subset \Pd\times
\Pd \times V$ is defined by the equation $\al(x)=0$. Therefore
$\phi^{-1}(W)$ is a hypersurface in the irreducible variety
$W\times X$ and any irreducible component  of $\phi^{-1}(W)$ has
the dimension $\dim X + \dim W - 1$ ($\phi^{-1}(W)$ is not empty
as $X$ meets any hyperplane and $\phi^{-1}(W)$ does not coincide
with $W\times X$ as $X$ is not contained in a hyperplane).

Thus the morphism $\phi$ is open. This proves Theorem 2.
\end{proof}

\begin{lemma} {\it Under the conditions of Theorem 2 the variety $M$ is
irredu\-cib\-le.}\end{lemma}

\begin{proof} It is sufficient to prove that the variety $M'\subseteq
V^*\times X$, $M' = \{(\al,x):\al(x)=0\}$, is irreducible.

(1) Consider the variety $L\subseteq V^*\times V$,
$L=\{(\al,v):\al(v)=0\}$, and the projection $\psi:L\to V$. The
morphism $\psi':L\setminus\psi^{-1}(0)\to V\setminus(0)$  is a
vector bundle. Indeed, fix a basis in $V$ and the dual basis in
$V^*$. Let $U_i=\{v\in V:x_i\ne 0\}$. Then $\psi'^{-1}(U_i)\cong
W\times U_i$, where $W$ is a vector space of dimension $n-1$ and
the isomorphism is defined by the following formula:
$$(a_1,\ldots,a_{i-1},a_{i+1},\ldots,a_n),(x_1,\ldots,x_n))\to$$
$$((a_1,\ldots,a_i=\frac{-1}{x_i}\sum_{j\ne i}
a_{j}x_j,\ldots,a_n),(x_1,\ldots,x_n)).$$

(2) Let $Z=M'\setminus\psi^{-1}(0)=\psi'^{-1}(X\setminus\{0\})$.
The map $\psi'':Z\to X\setminus\{0\}$ is also a vector bundle. The
variety $X$ is irreducible, consequently, $X\setminus\{0\}$ is
irreducible and $Z$ is irreducible.

(3) We have $M'=Z\sqcup Y$, where $Y$ is isomorphic to $V^*$ and
closed, $Z$ is irreducible. Let $M'=\cup M_i$, where $M_i$ are the
irreducible components of $M'$. It can be assumed that
$\overline{Z}\subseteq M_1$, $Y\subseteq M_2$. This yields that
there are no more then two irreducible components. If $Y\ne M_2$,
then $M_2=Y\cup(M_2\cap\overline{Z})$ and we have a contradiction
with the irreducibility of $M_2$. Hence $Y$ is an irreducible
component of $M'$. Further, $M'$ is a hypersurface in the
irreducible variety $X\times V^*$ and $\dim M_i=\dim M'$. At the
same time $\dim M'=n+\dim X-1>n$ and $\dim Y=n$. So we have a
contradiction. This means that $M'$ is irreducible. \end{proof}

Now we give some examples when the conditions of Theorem 2 do not
hold and $Ch(X)$ is not closed.

\begin{example} (1) $X$ is a line $\la v\ra\subset k^n, n\geq 2$ (here
$Ch(X)$ contains the open subset$\{(\la\al\ra,\la\be\ra)\ | \
\al(v)\ne~0\}$ of ${\bf P}(V^*)\times{\bf P}(V^*)$);

(2) $X$ is the subvariety in $k^3$ defined by the equations
$x_1=1$ and $x_2=x_3$ (indeed, $Ch(X)$ contains the subset
$\{(\la\al=a_1x_1+a_2x_2+a_3x_3\ra,\la\be\ra)\ |\ a_1\ne 0,
a_2+a_3=0\}$, the pair $(\la x_2-x_3\ra, \la x_1\ra)$ is contained
in its closure and is not contained in $Ch(X)$);

(3) $X\subset V=k^2$ is defined by the equation $x_1x_2=0$ (here
$Ch(X)$ contains the open subset
$\{(\la\al=a_1x_1+a_2x_2\ra,\la\be\ra)\ | \ a_1\ne 0, a_2\ne 0\}$
of ${\bf P}(V^*)\times{\bf P}(V^*)$).\end{example}

Theorem 2 implies the similar theorem for a subset in a projective
space (the definition of the characteristic variety in the
projective case is analogous).

\medskip
\noindent {\bf Theorem 2'}. {\it Suppose that an algebraic
subvariety $Y\subset \PProj$ is irreducible and is not contained
in a hyperplane. Then $Ch(Y)$ and $Chw(Y)$ are closed.}

\begin{proof} Let $X\subset V$ be the cone corresponding to $Y$.
Note that $Ch(X)=Ch(Y)$. Applying  Theorem 2 to $X$, we conclude
that $Ch(X)$ is closed.\end{proof}

\section{The case of a $T$-invariant subvariety}

Let $T$ be an algebraic torus linearly acting on a vector space
$V$ and $X$ be a $T$-invariant subset in $V$. Then $Ch(X)$ is a
$T$-invariant subset in $\Pd\times\Pd$.

\begin{proposition} {\it If $X$ has not the separation property and $Ch(X)$ is
closed, then there exists a pair $(\la\al\ra,\la\be\ra)\in Ch(X)$
such that $\al, \be$ are $T$-semiinvariant and linearly
independent. }\end{proposition}

\begin{proof} Fix a $T$-semiinvariant basis $\{x_1,\ldots,x_n\}$ in $V^*$.
Let $\lambda_i$ be the weight of $x_i$.

There exists a pair $(\la\al\ra,\la\be\ra)\in Ch(X)$ such that
$\la\al\ra\ne \la\be\ra$. Consider the action
$T:\overline{T(\la\al\ra,\la\be\ra)}$. For an action of a torus on
a projective variety there exists a fixed point. Denote it by
$(\la\al'\ra,\la\be'\ra)$. If $\la\al'\ra\ne\la\be'\ra$, then
there is nothing to prove. Now let $\la\al'\ra=\la\be'\ra$.

There exists a one-parameter subgroup $\ga: k^* \to T$ such that
$$\lim_{s\to 0}\ga(s)(\la\al\ra,\la\be\ra)=(\la\al'\ra,\la\al'\ra)$$
(see \cite[Sec. 2.3]{Ful}).

Let $(\al_1, \ldots, \al_n)$ be the coordinates of $\al$, and
$(\be_1, \ldots, \be_n)$ be the coordinates of $\be$. We may
assume that $\al' = (\al'_1,\ldots \al'_n)$, where $\al'_i=\al_i$
if $\la\ga,\lambda_i\ra=\min_{j\ :\ \al_j\ne
0}\la\ga,\lambda_j\ra$, and $\al'_i=0$ in another case, and the
formula for $\be'$ is analogous (here $\la\cdot\ ,\ \cdot\ra$ is
the natural pairing between the lattice of one-parameter subgroups
and the lattice of characters ). Let $\be'=c\al'$.

We have $(\la\al\ra, \la\be-c\al\ra)\in Ch(X)$ and $$\lim_{s\to
0}\ga(s)(\la\al\ra,\la\be-c\al\ra)=(\la\al'\ra,\la\be''\ra)\in
Ch(X).$$ Note that the supports of the vectors $\al',\be''$ do not
intersect. (Here the support of a vector $v$ is the set of basis
vectors along which $v$ has  non-zero coordinates; the support of
$\la v \ra$ is the support of $v$.) Thus we obtain
$\overline{T(\la\al'\ra,\la\be''\ra)}\cap D=\emptyset$ and the
desired point is a $T$-fixed point in
$\overline{T(\la\al'\ra,\la\be''\ra)}$. \end{proof}

\begin{corollary} Let $T:V$ and all $\lambda_i$ be different. Let $\{\e_i\}$
be a $T$-semiinvariant basis in $V$, $\{x_i\}$ be the dual basis
in $V^*$, and $X\subset V$ be a closed irreducible affine
$T$-invariant subvariety. Then $X$ does not satisfy (SP) $\iff$
$X$ is contained in a hypersurface of the form $x_i^n + x_jF=0$
for some $i\ne j$, $F\in k[V]$ and $n\in \N$. \end{corollary}

\begin{proof} We have to prove that the separation property fails
on a pair of $T$-semiinvariant linear functions.

(1) If $X$ is contained in a homogeneous hyperplane, then $X$ is
contai\-ned in a $T$-invariant homogeneous hyperplane and (SP)
fails on a $T$-semiinvariant pair.

(2) If $X$ does not meet a homogeneous hyperplane, then there
exist $f, g\in k[X]$, such that $fg=1$. If the functions $f$ and
$g$ are not $T$-semiinvariant, then we consider their weight
decompositions. Since $k[X]$ has no zero divisors, it follows that
after the multiplication in the left side of the equality we
obtain the sum of weight functions with different weights. There
exists a one-parameter subgroup having different pairing with all
weights from the weight decompositions of $f$ and $g$. This
one-parameter subgroup defines the order on weight functions. The
products of the highest and the lowest terms can not be cancelled,
so we have a contradiction. Thus $X$ does not meet a $T$-invariant
homo\-ge\-neous hyperplane.

(3) If $\dim X = 1$, then either $X$ is a curve of $T$-fixed
points or $X$ is a closure of an orbit of a one-dimensional torus.
In the first case $X$ is contained in the subspace defined by the
equations $x_i=0$, where $\lambda_i\ne 0$. In the second case (SP)
fails on any pair of coordinate functions $x_i, x_j$ such that
either $\lambda_i=0$ or $\lambda_i, \lambda_j\ne 0$.

(4) If $X$ is not contained in a homogeneous  hyperplane, meets
any homogeneous hyperplane and $\dim X>1$, then, by Theorem 2,
$Ch(X)$ is closed and, by Proposition 3, (SP) fails for $X$ on a
pair of $T$-semiinvariant functions.\end{proof}

\section{Application to binary forms}

Let ${\rm char}\ k=0.$  Consider the vector space $k[x,y]_n$ of
binary forms of degree $n$, where $SL_2$ acts by the natural way
and $k^*$ acts by homotheties.

K. Baur proved the following theorem in \cite{Baur}.

\medskip
\noindent {\bf Theorem \cite[Th. 3.4]{Baur}.} {\it Let $f\in
k[x,y]_n$. Then the orbit $O_f=SL_2 f$ satisfies (SP) if and only
if $f$ has a linear factor of multiplicity one.}

\medskip
We give the proof of this theorem here. Corollary 1 allows to
simplify it (see Proposition 4 below).

The proof consists of some lemmas.

\begin{lemma} {\it Suppose that the orbit $O_f$ satisfies (SP). Then $f$ has a
linear factor of multiplicity one.}\end{lemma}

\begin{proof} Suppose that any linear factor of $f$ has the
multiplicity no less then two. Note that this property holds for
any $h(x,y)=h_0 x^n+h_1 x^{n-1}y+\ldots+h_n y^n\in O_f$. Then (SP)
fails for $O_f$ on the pair of linear functions $\al(h)=h_0,
\be(h)=h_1$. Indeed, if $h_0=0$ for some $h\in O_f$, then $h$ is
divisible by $y$. This implies that $h$ is divisible by $y^2$ and
$h_1=0$.\end{proof}

\begin{lemma}{\it The orbit $O_f$ satisfies (SP) if and only if the closure
$\overline{k^*O_f}$ satisfies (SP).}\end{lemma}

\begin{proof} If $O_f$ has (SP), then $\overline{k^*O_f}$ has
(SP) (as $O_f\subset \overline{k^*O_f}$).

Suppose that $O_f$ has not (SP), namely there exist homogeneous
hyperplanes $H\ne H'$ such that $H\cap O_f\subset H'$. Then $H\cap
tO_f=t(H\cap O_f)\subset H'$ for any $t\in k^*$ and $H\cap
k^*O_f\subset H'$. By Lemma 4 below, it follows that
$H\cap\overline{k^*O_f}=\overline{H\cap k^*O_f}\subset H'$ and
$\overline{k^*O_f}$ has not (SP). \end{proof}

\begin{lemma} \cite[Lemma 3(c)]{KrW}\ {\it Let $G:V$ be an
irreducible representation of a connected algebraic group $G$,
$H\subset V$ be a homogeneous hyperplane, and $X\subset V$ be a
constructive $G$-invariant subset. Then $\overline{H\cap
X}=H\cap\overline{X}$.}\end{lemma}

\begin{proposition} {\it The orbit $O_{xy^{n-1}}$ satisfies (SP).}\end{proposition}

\begin{proof} It is sufficient to prove that
$$X=\overline{k^*SL_2 xy^{n-1}}=\overline{GL_2 xy^{n-1}}=\{(ax+by)(cx+dy)^{n-1}\ |\ a,b,c,d\in k\}$$
has the separation property (see Lemma 3). For a point
$h=(ax+by)(cx+dy)^{n-1}\in X$ we have
$$h_m=aC_{n-1}^{m-1}c^{m-1}d^{n-m}+bC_{n-1}^{m}c^{m}d^{n-m-1}=$$
$$=\frac{(n-1)!c^{m-1}d^{n-m-1}(mad+(n-m)bc)}{m!(n-m)!}.$$

If the separation property does not hold, then $X$ is contained in
a hypersurface of the form $h_i^k + h_jF=0$,  where $F$ is a
polynomial in the variables $h_m$ (see Corollary 1). Putting in
this equation coordinates of points from $X$, we get a polynomial
in the variables $a,b,c,d$ identically equal to zero. This implies
that $h_i^k$ is divisible by $h_j$ (as a polynomial in $a,b,c,d$).
It is easy to see that this can not be true for any
$i,j,k$.\end{proof}

The following lemma combining with Lemma 3 completes the proof of
the Theorem.

\begin{lemma} {\it If $f\in k[x,y]_n$ has a linear factor of
multiplicity one, then $\overline{k^*O_f}\supseteq
O_{xy^{n-1}}$.}\end{lemma}

\begin{proof} We may assume that $f$ has the form $x(y+a_2x)\ldots$
$(y+a_nx)$. Let $\gamma(t)$ be the diagonal matrix $\left(
\begin{array}{cc}
t& 0 \\
0 & t^{-1}
\end{array}
\right)\in SL_2$. We have
$$t^{n-2}\gamma(t)f=x(y+t^2a_2x)\ldots(y+t^2a_nx)\rightarrow xy^{n-1}
 \mbox{ as } t\rightarrow 0.$$\end{proof}

\begin{remark}\ \ Lemma \ 3 \ \ implies \ \ the \ \ similar \ \ result\ \  for\ \ the\ \ action\ \ $SL_2 :
\mathbb{P}(k[x,y]_n)$. Let $f\in \mathbb{P}(k[x,y]_n)$. Then the
orbit $O_f=SL_2 f$ satisfies (SP) if and only if $f$ has a linear
factor of multiplicity one.\end{remark}

\section{The closure of a toric orbit in a vector space}

Let $T$ be an algebraic torus and $\Lambda$ be the lattice of its
characters. Consider a linear action $T:V$, where
$t(x_1,\ldots,x_n)=(\chi_1(t)x_1,\ldots,\chi_n(t)x_n)$. Let
$\Sigma$  be the semigroup in $\Lambda$ generated by the
characters $\chi_1,\ldots,\chi_n$ and $K$ be the cone in
$\Lambda\otimes_{\Z}\Q$ generated by $\Sigma$.

We are interested in the question if $X=\overline{Tv}$ has the
separation property. If the vector $v$ has a zero coordinate, then
$X$ is contained in a homogeneous hyperplane and has not (SP). If
any coordinate of $v$ is non-zero, then, up to the proportionality
of the basis vectors, we may assume that $v=(1,\ldots,1)$. It can
also be assumed that the kernel of inefficiency of the action
$T:V$ is finite (i.e., $\dim X=\dim T$, or the cone $K$ is
bodily).

\begin{lemma} \cite[Lemma 4.1]{St}{\it Let $X=\overline{T(1,\ldots,1)}$.
The ideal $I(X)$ is generated (as a vector space) by all binomials
of the form $x_1^{a_1}\cdot\ldots\cdot
x_n^{a_n}-x_1^{b_1}\cdot\ldots\cdot x_n^{b_n}$, where $\sum
a_i\chi_i = \sum b_i\chi_i$, $a_i, b_i \in \Z_+$.}
\end{lemma}

Now we give an algorithm for a construction of a finite system of
binomials defining $X$. Denote by $W$ the sublattice in $\Z^n$
defined by the system of equations $\sum c_i\chi_i=0$, where
$c_i\in \Z$. To any point $c=(c_1,\ldots,c_n)$ of this sublattice
we put in correspondence the binomial $x_1^{a_1}\cdot\ldots\cdot
x_n^{a_n}-x_1^{b_1}\cdot\ldots\cdot x_n^{b_n}$, where $a_i=c_i$,
$b_i=0$ if $c_i>0$, and $b_i=|c_i|$, $a_i=0$ otherwise. Such
binomials vanish on $X$. Consider an octant $D_I=\{c\in \Z^n \ | \
c_i\geq 0$ if $i\in I$, $c_i\leq 0$ if $ i\not \in I\}$. Then
$W\cap D_I$ is a finitely generated semigroup. Let us construct a
system of generators of $W\cap D_I$ (for example, if
$d_1,\ldots,d_k\in W\cap D_I$ is a \ system \ of \ generators\ of
the \ cone \ $\Q_+(W\cap D_I)$,\ then \ the \ set \ $\{\sum s_i
d_i \ | \ 0\leq s_i\leq 1\}\cap D_I$ generates $W\cap D_I$). Now
we consider the union of such systems of generators of the
semigroups $W\cap D_I$ over all octants and prove that the
corresponding set of binomials generates the ideal $I(X)$. Indeed,
the element $x^{I_1+J_1}-x^{I_2+J_2}$ belongs to the ideal
generated by $x^{I_1}-x^{I_2}$ and $x^{J_1}-x^{J_2}$:
$$x^{I_1+J_1}-x^{I_2+J_2}=x^{J_1}(x^{I_1}-x^{I_2})+x^{I_2}(x^{J_1}-x^{J_2})$$
(here $I_1, I_2, J_1, J_2$ are ordered sets of $n$ non-negative
integers and $x^{J}=x_1^{c_1}\ldots x_n^{c_n}$, where
$J=(c_1,\ldots, c_n)$, $c_i\in \mathbb{Z}_+$). Hence the binomials
corresponding to the constructed system of generators of $W$ as a
semigroup generate the ideal $I(X)$.

\begin{theorem} {\it The orbit closure $X=\overline{T(1,\ldots,1)}$
satisfies (SP) if and only if the following conditions hold:

\noindent (1) the cone $K$ is strictly convex;

\noindent (2) $\forall i \quad \Q_+\chi_i$ is an edge of $K$;

\noindent (3) $\Q_+\chi_i\ne\Q_+\chi_j$ if $i\ne j$.}
\end{theorem}

\begin{remark} It is easy to see that these conditions are equivalent to
the following:

\noindent (a) $\chi_i\not\in
<\chi_1,\ldots,\chi_{i-1},\chi_{i+1},\ldots,\chi_n>_{\Q_+}$ for
any $i$;

\noindent (b) $-\chi_i\not\in
<\chi_1,\ldots,\chi_{i-1},\chi_{i+1},\ldots,\chi_n>_{\Q_+}$ for
any $i$.
\end{remark}

\begin{proof} (1) If $K$ contains a line, then there exist
$c_1,\ldots,c_n\in \Z_+$ such that $c_i\ne 0$ for some $i$ and
$\sum c_i\chi_i=0$. Therefore $X$ is contained in a hypersurface
$x_1^{c_1}\cdot\ldots\cdot x_n^{c_n}=1$ and does not meet the
hyperplane $x_j=0$, where $c_j\ne 0$. Consequently, $X$ has not
(SP).

(2) If $\mathbb{Q}_+\chi_i$ is not an edge of $K$, then there
exist $c_1,\ldots,c_n\in \Z_+$ such that $c_i\chi_i = \sum_{j\ne
i}c_j\chi_j$, $c_i\in\N$. If $\chi_i=0$, then $X$ does not meet
the hyperplane $x_i=0$ and has not (SP). If $\chi_i\ne 0$, then
there exists $c_j\ne 0$, $j\ne i$. This means that $x\in X, x_j=0$
implies $x_i=0$ and (SP) does not hold. The case
$\Q_+\chi_i=\Q_+\chi_j$ can be proved by the same arguments.

(3) Now we prove the converse implication. Consider the case $\dim
T=1$. The set $\{\chi_1,\ldots,\chi_n\}$ satisfies the conditions
of the theorem only if $\dim V=1$. In this case $X=V$ if
$\chi_1\ne 0$ and $X=\{1\}$ if $\chi_1 = 0$ and $X$ has (SP).

Now let $\dim X>1$. Note that $X$ is not contained in a hyperplane
(otherwise there exist $i\ne j$ such that $\chi_i=\chi_j$) and
meets any hyperplane (since $K$ is strictly convex and $\chi_i\ne
0$ for any $i$, it follows that $0\in X$). Hence, by Theorem 2,
$Ch(X)$ is closed. If $X$ has not (SP), then, by Corollary 1, $X$
is contained in the hypersurface $x_i^l + x_jF$ for some $i\ne j,
F\in k[V], l\in \mathbb{N}$. By Lemma 6, it follows that $x_i^l
-x_1^{c_1}\ldots x_n^{c_n}\in I(X)$ for some $c_1,\ldots,c_n\in
\mathbb{Z}_+$, $c_j
> 0$, i.e., $l\chi_i=\sum c_m\chi_m$. If $l\leq c_i$, then since
$\chi_j\ne 0$, we have a contradiction with condition (1).
Otherwise we have a contradiction with $\chi_i\not\in
<\chi_1,\ldots,\chi_{i-1},\chi_{i+1},\ldots,\chi_n>_{\Q_+}$.
\end{proof}

\begin{corollary}\ Let \ the orbit\ closure\
$X=\overline{T(1,\ldots,1)}$ \ be \ defined \ by \ $p_1=0, \ldots,
p_m=0$, where $p_i$ are binomials, and $S_i$ be the hypersurface
defined by $p_i=0$. Then the following conditions are equivalent:

\noindent (1) $X$ does not satisfy  (SP);

\noindent (2) there exists $i$ such that $S_i$ does not satisfy
(SP);

\noindent (3) there exists $i$ such that $p_i$ has either the form
 $x_1^{c_1}\ldots x_n^{c_n}-1$,
where $c_k\geq 0$ and  some $c_j>0$, or the form
$x_i^{c_i}-x_1^{c_1}\ldots
x_{i-1}^{c_{i-1}}x_{i+1}^{c_{i+1}}\ldots x_n^{c_n}$, where
$c_k\geq 0, c_i>0$.
\end{corollary}

\begin{proof} The implications $(3)\Rightarrow (2)$ and $(2)\Rightarrow
(1)$ are obvious.

Suppose that among $p_i$ there are no equations of such form. We
shall prove that $X$ has  (SP) applying Theorem 3.

(a) The cone $K$ is strictly convex and $\chi_i\ne 0$ for any $i$
since $0\in X$.

(b) Suppose that $\Q_+\chi_i$ is not an edge of $K$ or
$\Q_+\chi_i=\Q_+\chi_j$. This means that some equation of the form
$x_i^{c_i}=x_1^{c_1}\ldots
x_{i-1}^{c_{i-1}}x_{i+1}^{c_{i+1}}\ldots x_n^{c_n}$, where
$c_k\geq 0, c_i>0$ and some $c_j>0, j\ne i$, vanishes on $X$. On
the other hand, $(0,\ldots,x_i=1,\ldots,0)\in X$ and we have a
contradiction. By Theorem 3, it follows that $X$ has  (SP) and the
implication $(1)\Rightarrow (3)$ is proved. \end{proof}

\begin{example} Let $T=(k^*)^3$. Consider the 5-dimensional representation
of $T$ with the characters $\chi_1=(1,0,0),$ $\chi_2=(1,1,0),$
$\chi_3=(0,1,2),$ $\chi_4=(0,2,1),$ $\chi_5=(1,0,1).$ (On the
picture below we draw the  corresponding cone $K$.) Then
$X=\overline{T(1,1,1,1,1)}$ can be defined by the equations
\begin{equation*}
\left\{%
\begin{array}{c}
  x_1^3x_3=x_2x_5^2, \\
  x_1^3x_4=x_2^2x_5, \\
  x_2x_3=x_4x_5.
\end{array}
\right.
\end{equation*}

The characters of this representation satisfy the conditions of
Theorem 3 (and the equations do not satisfy the conditions of
Corollary 2). Thus $X$ has the separation property.
\end{example}

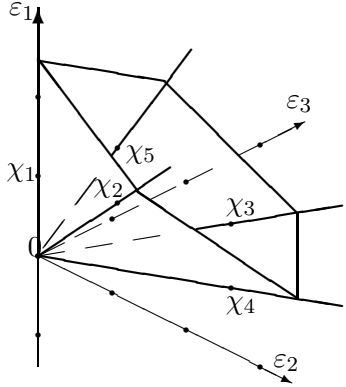
\begin{picture}(140,150)
\put(21,7){\line(0,1){42}} 
\multiput(21,19)(0,30){4}{\circle*{2}} \put(17,49){$0$}
\put(9,79){$\chi_1$} \put(10,140){$\varepsilon_1$}



\put(21,49){\vector(2,-1){96}} 
\multiput(21,49)(28,-14){4}{\circle*{2}}
\put(110,7){$\varepsilon_2$}

\put(97,87){\vector(2,1){25}} 
\multiput(21,49)(28,14){4}{\circle*{2}}
\put(115,105){$\varepsilon_3$}
\multiput(21,49)(14,7){5}{\line(2,1){10}}


\thicklines

\put(21,49){\vector(0,1){94}}

\put(21,49){\line(6,-1){115}}

\put(81,59){\line(6,1){55}}

\put(21,49){\line(3,2){50}}

\put(49,87){\line(3,4){30}}

\thinlines

\multiput(23,52)(12,16){2}{\line(3,4){8}}
\multiput(21,49)(18,3){3}{\line(6,1){10}}

\put(94,37){\circle*{2}} \put(92,29){$\chi_4$}
\put(94,61){\circle*{2}} \put(92,66){$\chi_3$}
\put(51,90){\circle*{2}} \put(54,86){$\chi_5$}
\put(51,69){\circle*{2}} \put(42,73){$\chi_2$}


\thicklines

\put(119,33){\line(0,1){32}}

\put(119,33){\line(-3,2){61}}

\put(59,73){\line(-3,4){38}}

\put(21,123){\line(6,-1){48}}

\put(69,115){\line(1,-1){50}}

\end{picture}

\begin{remark} Note that an orbit of a torus has not (SP) ($x_i\ne 0$
on $T(1,\ldots,1)$). \end{remark}

\begin{remark} We say that a variety $Y\subset V$ is {\it binomial} if $Y$ can be
defined by a system of binomials. In particular, the closure of an
orbit of a torus is binomial. Moreover, it is easy to see that a
binomial variety is the closure of some toric orbit if and only if
it is irreducible.  The following example shows that the statement
of Corollary 2 does not hold for an arbitrary binomial
variety.\end{remark}

\begin{example} Consider the binomial variety $Y\subset k^4$ defined by the
equations
\begin{equation*}
\left\{%
\begin{array}{c}
  x_1x_2=x_3x_4, \\
  x_1x_3=x_2x_4, \\
  x_1x_4=x_2x_3. \\
\end{array}
\right.
\end{equation*}

Any hypersurface $x_ix_j=x_lx_m$ has (SP). It is easy to check
that
 (SP) (and even (WSP)) for $Y$ fails on the functions
$x_1-ax_2$ and $x_1-bx_2$, where $a,b\ne\pm 1, a\ne b$ (since $Y$
is an union of four coordinate axes and the lines $x_1=\pm x_2=\pm
x_3=\pm x_4$ with an even number of minuses).\end{example}

\begin{theorem} {\it The orbit closure
$X=\overline{T(1,\ldots,1)}$ satisfies (WSP) if and only if the
following conditions hold:

\noindent (1) the cone $K$ is strictly convex;

\noindent (2) there are no more then one characters $\chi_i$ in
the interior of any face of $K$ (in particular,
$\Q_+\chi_i\ne\Q_+\chi_j$ for $i\ne j$).}\end{theorem}

\begin{proof} (1) If $K$ contains a line, then there exist $c_1,\ldots,$
$c_n\in \Z_+$ such that $c_j, c_k\ne0$ for some $j\ne k$ and $\sum
c_i\chi_i=0$. Consequently, $X$ is contained in the hypersurface
$x_1^{c_1}\cdot\ldots\cdot x_n^{c_n}=1$ and does not meet the
hyperplanes $x_j=0$ and $x_k=0$. Therefore $X$ has not (WSP).

(2) Let condition (2) do not hold, i.e, there are two characters
(let us assume that they are $\chi_1, \chi_2$) in the interior of
the face of $K$ generated by $\{\chi_i \ | \ i\in I\}$. Then there
exist $c_i, d_i>0$ such that $c_1\chi_1=\sum_{i\in I} c_i\chi_i$
and $d_2\chi_2=\sum_{i\in I} d_i\chi_i$. We have
$$X\cap H_{x_1}=\cup_{i\in I}X\cap H_{x_i}=X\cap H_{x_2}$$ and $X$ has
not (WSP).

(3) Let us prove the converse implication. Consider the case $\dim
T=1$. The cone $K$ satisfies the conditions of the theorem if and
only if either $\dim V=1$ or $\dim V=2$ and $t(x_1,x_1)=(x_1,t^m
x_2)$, where $m\ne 0$. In the first case $X$ has (SP) and it is
easy to see that in the second case $X$ has (WSP).

Now let $\dim T\geq 2$. Assume the converse. Let the cone $K$
satisfy the conditions of Theorem 4 and the weak separation
property fails for $X$ on  linearly independent functions
$\al=a_1x_1+\ldots+a_nx_n, \be=b_1x_1+\ldots+b_nx_n$.

If $M$ is a cone generated by some edges of $K$ and $\delta$ is a
linear function, then denote by $\delta_M$ the restriction of
$\delta$ to the subspace in $V$ defined by the equations $x_i=0$,
where $\chi_i\not \in M$.

\begin{remark} Let $L$ be a proper face of $K$ and $W \subset V$ be the
subspace defined by the equations $x_i=0$, where $\chi_i\not \in
L$. Note that $X\cap W$ is the closure of the $T$-orbit of the
point with the coordinates $x_i=1$ for $\chi_i\in L$ and $x_i=0$
for $\chi_i\not \in L$. Since the characters of the representation
$T:V$ satisfy the conditions of the theorem, it follows that the
characters of the representation $T:W$ satisfy the same
conditions. By inductions over $\dim V$, we may assume that $X\cap
W\subset W$ has (WSP). Since $H_{\al_L}\cap X\cap W=H_{\be_L}\cap
X\cap W$, it follows that $\al_L$ and $\be_L$ are linearly
dependent.
\end{remark}

\begin{remark} If the linear functions $\tilde{\al}(x)=a_1x_1+\ldots
+a_{i-1}x_{i-1}+a_{i+1}x_{i+1}+\ldots+a_nx_n$ and $\tilde{\be}(x)=
b_1x_1+\ldots +b_{i-1}x_{i-1}+b_{i+1}x_{i+1}+\ldots+b_nx_n$ are
proportional for some $i$, then for $x\in X$ we have $\al(x)=0
\iff \be(x)=0 \iff \tilde{\al}(x)=0$ and $x_i=0$. There exists a
one-parameter subgroup $\gamma:k^*\to T$ such that $\la \gamma,
\chi_r\ra\ne \la \gamma, \chi_j\ra$ for $r\ne j$. Hence there
exists a non-zero root of the equation $a_1t^{\la \gamma,
\chi_1\ra}+\ldots +a_nt^{\la \gamma, \chi_n\ra}=0$. Thus we have a
contradiction.
\end{remark}

The first case. Suppose that $\chi_i\ne 0$ for any $i$.
\medskip

\noindent {\bf Step 1.} Suppose that $\al(x)=a_ix_i+a_jx_j$ and
$\be(x)=b_ix_i+b_jx_j$. Note that $a_ix_i+a_jx_j=0$ and
$b_ix_i+b_jx_j=0$ if and only if $x_i=0,x_j=0$. Since $\chi_i\ne
\chi_j$, it follows that there exists $x\in T(1,\ldots,1)$ such
that $\al(x)=0$ and we have a contradiction.

\medskip \noindent {\bf Step 2.} Suppose that $\al(x)=a_ix_i+a_jx_j+a_mx_m$
and $\be(x)=b_ix_i+b_jx_j+b_mx_m$, where $\chi_m$ belongs to the
interior of $K$.

\noindent 1. If $a_i=b_i=0$ or $a_j=b_j=0$, then we have a
contradiction (see Step~1).

\noindent 2. If $a_i=0, b_i\ne 0$, then there exists a
one-parameter subgroup $\gamma: k^*\to T$ such that
$\gamma(s)\al\to \al'$, $\gamma(s)\be\to \be'$, where the supports
of $\al'$ and $\be'$ do not intersect. Since $Chw(X)$ is closed,
it follows that $(\la \al'\ra,\la \be'\ra)\in Chw(X)$ and
$\overline{T(\la \al'\ra, \la \be'\ra)}\cap D=\varnothing$. There
exists a $T$-stable point $(\la x_l \ra,\la x_r\ra)\in
\overline{T(\la \al'\ra, \la \be'\ra)}$. Consequently, $X$ is
contained in a hypersurface of the form $x_l^{c}+x_r^{d}+x_lx_rF$
for some $l\ne r, c,d\in \mathbb{N}, F\in k[X]$. Lemma 5 implies
that the binomials $x_l^{c}-x_1^{c_1}\ldots x_n^{c_n}$ and
$x_r^{d}-x_1^{d_1}\ldots x_n^{d_n}$ belong to $I(X)$ for some
$c_l, d_l$ such that $c_r, d_l>0$ and $c\chi_l=\sum c_i\chi_i$,
$d\chi_r=\sum d_i\chi_i$. If $c\leq c_l$ or $d\leq d_r$, then this
contradicts condition (1). Otherwise $\chi_l$ and $\chi_r$ belong
to the interior of the same face of $K$ and we also have a
contradiction. The cases when $a_i\ne 0, b_i=0$ or $a_j=0, b_j\ne
0$ or  $a_j\ne 0, b_j=0$ can be considered by the same way.

\noindent 3. Now let $a_i,b_i,a_j,b_j\ne 0$. It can be assumed
that $a_i=b_i=1$. There exists a one-parameter subgroup $\gamma:
k^*\to T$ such that $\la \gamma, \chi_i \ra = n_1 > \la \gamma,
\chi_m \ra = n_2 > \la \gamma, \chi_j \ra = 0$.
The equations $$s^{n_1}+a_ms^{n_2}+a_j=0 \ \ \ \ \ \ \ \ \ \ (*)$$
and $$\ s^{n_1}+b_ms^{n_2}+b_j=0\ \ \ \ \ \ \ \ \ \ (**)$$ have
the same roots up to multiplicities. If all  roots of one equation
have the multiplicity one, then all roots of another equation have
the multiplicity one and the equations are proportional. A
multiple root $s_0$ of $(*)$ is a root of its derivative and
satisfies the equation
$$s_0^{n_1-n_2}=-\frac{a_mn_2}{n_1}.$$ (If $n_1$ is divisible by
$char\ k$, then either $a_m=b_m=0$ and this is Step 1 or $n_2$ is
divisible by $char \ k$. In the second case we replace $n_1$ by
$\frac{n_1}{char\ k}$, $n_2$ by $\frac{n_2}{char\ k}$, and
$a_j,a_m,b_j,b_m$ by the roots of degree $char\ k$ and  obtain
similar equations.) Putting obtained formul\ae\ in  $(*)$, we get
$$s_0^{n_2}(-\frac{a_mn_2}{n_1})+a_ms_0^{n_2}+a_j=0.$$ This
implies $$s_0^{n_2}=-\frac{a_jn_1}{a_m(n_1-n_2)}$$ (if $a_m=0$ or
$n_1-n_2$ is divisible by $char\ k$, then $n_1$ is divisible by
$char\ k$ and we repeat the previous arguments). Also, we have
$$s_0^{n_1}=\frac{a_jn_2}{n_1-n_2}.$$ Note that $s_0$ is a root of
$(**)$. Putting obtained formul\ae \ in $(**)$, we get
$$a_ma_jn_2-b_ma_jn_1+b_ja_m(n_1-n_2)=0.$$ Since $(**)$
also has a multiple root, then, by the same arguments, we have the
symmetric equation
$$b_mb_jn_2-a_mb_jn_1+a_jb_m(n_1-n_2)=0.$$ Summing the obtained
equations, we get $(a_j-b_j)(a_m-b_m)=0$. Thus we have a
contradiction with Remark 11.

\medskip

\noindent {\bf Step 3.} Consider the case of arbitrary $\al,\be$.
In the case $\dim T = 2$ we have a contradiction with Step 2.

Now let $\dim T>2$. We say that $\tilde{K}$ is a {\it subcone} of
$K$ if $\tilde{K}$ is generated by a finite number of vectors from
$K$. Note that the interior of a subcone is contained in the
interior of one of the faces of $K$. Consider the subcone $K'$ of
$K$ generated by the characters $\chi_i$ such that $a_i\ne 0$ or
$b_i\ne 0$. Since $\al_{K'}$ and $\be_{K'}$ are not proportional,
it follows that there exists a face $L_1$ of codimension 1 in $K'$
such that $\al_{L_1}$ and $\be_{L_1}$ are not proportional. If
$\chi_m$ is contained in the interior of $L_1$, then $\chi_m$ is
contained in the interior of $K$. Indeed, assume the converse. Let
$\chi_i$ be an interior point of a proper face of $K$ and $\chi_i$
belong to the interior of $L_1$. Then $L_1$ is contained in this
face of $K$. By Remark 10, the linear functions $\al_{L_1}$ and
$\be_{L_1}$ are proportional and we have a contradiction. For the
same reason, there exists a face $L_2$ of codimension one in $L_1$
such that $\al_{L_2}$ and $\be_{L_2}$ are   not proportional. Thus
we can find a 2-dimensional face $L_r$ with this property. There
exists a one-parameter subgroup $\gamma:k^*\to T$ such that
$\gamma(s)\al\to \al_{L_r}, \gamma(s)\be\to \be_{L_r}$ as $s\to
0$. Since $Chw(X)$ is closed, it follows that $(\la
\al_{L_r}\ra,\la \be_{L_r}\ra)\in Chw(X)$ and we can apply  Step
2.

\medskip The second case. Suppose that $\chi_i=0$ for some $i$
(we may assume that $i=1$). If $a_1=b_1=0$, then consider the
image $X'\subset W$ of $X$ under the projection along the first
basis vector (here $W$ is a subspace defined by the equation
$x_1=0$). The characters of the representation $T:W$ are non-zero
and satisfy the conditions of the theorem. The first case implies
that $X'\subset W$ has  (WSP). But  (WSP)  fails for $X'\subset W$
on the restriction on $W$ of the functions $\al,\be$. This is a
contradiction.

Let $a_1$ or $b_1$ be not equal to zero. By Remark 10, it follows
that the vector  $(a_i, b_i)$ is proportional to the vector $(a_1,
b_1)$ for $i\ne m$. Thus the linear functions $\al'=a_1x_1+\ldots
+a_{m-1}x_{m-1}+a_{m+1}x_{m+1}+\ldots+a_nx_n$ and $\be'=
b_1x_1+\ldots +b_{m-1}x_{m-1}+b_{m+1}x_{m+1}+\ldots+b_nx_n$ are
proportional. This contradicts to Remark 11. \end{proof}

\begin{theorem} {\it Suppose that the orbit closure
$X=\overline{T(1,\ldots,1)}\subset V$ is a cone; then $X$
satisfies (SSP) if and only if $X=V$ (i.e., the weights
$\chi_1,\ldots,\chi_n$ are linearly independent).} \end{theorem}

\begin{proof} Let $X\ne V$ and $x_1^{a_1}\ldots
x_n^{a_n}-x_1^{b_1}\ldots x_n^{b_n}\in I(X)$, where there exists
$i$ with $a_i\ne b_i$. We may assume that  $a_i=0$ or $b_i=0$ for
any $i$. Let $a_1>0$. Since $X$ is a cone, it follows that there
exists $b_i>0$ and $$X\cap H_{x_1}=\cup_{i:\ b_i>0} (X\cap V_i),$$
where\ \ $V_i = H_{x_1}\cap H_{x_i}$. \ This \ \ implies\ \ that\
there \ exists\  $i$ \ such \ that\ $\dim X\cap V_i = \dim X\cap
H_{x_1}$ and $X\cap V_i$ has the codimension $\leq 1$ in $X$.
\end{proof}

\begin{remark} The proof of Theorem 5 is true for any cone
which is contained in a binomial hypersurface, in particular, for
binomial cones.\end{remark}

\noindent {\bf Question.} Let $X$ be a closed irreducible
$T$-invariant subvariety in a $T$-module $V$ such that $X$ has not
(WSP) (resp. (SSP)). Is it true that  (WSP) (resp. (SSP)) for $X$
fails on a pair of $T$-semiinvariant linear functions?

\section{The closure of a toric orbit in a projective space}

Let $T:\PProj$,
$t(x_1:\ldots:x_n)=(\chi_1(t)x_1:\ldots:\chi_n(t)x_n)$. We are
interested in the question if $Y=\overline{Tw}$ ($w\in \PProj$)
has the separation properties. As in the previous section we may
assume that $w=(1:\ldots:1)$. It also can be assumed that the
kernel of inefficiency of the action $T:\PProj$ is finite, i.e.,
$\dim X=\dim T$.

\medskip
\noindent {\bf Theorem 3'.} {\it The orbit closure
$X=\overline{T(1:\ldots:1)}$ satisfies (SP) if and only if the
following conditions hold:

\noindent (1) the point $\chi_i$ is a vertex of the convex hull
${conv}\{\chi_1,\ldots,\chi_n\}$ for any $i$ ;

\noindent (2) $\chi_i\ne \chi_j$ for $i\ne j$.}

\begin{proof} Consider the action $T\times k^*:V$, where
$$(t,s)(x_1,\ldots,x_n) = (s\chi_1(t)x_1,\ldots,s\chi_n(t)x_n).$$
The weights of this representation are $\chi_i'=\chi_i+\lambda$,
where $\lambda(t,s)=s$. Then $$T(1:\ldots:1) =
{\mathbb{P}}(T\times k^* (1,\ldots,1)),$$ and
$X=\overline{T(1:\ldots:1)}={\mathbb{P}}\overline{(T\times
k^*(1,\ldots,1))}$. We shall apply Theorem 4 to $\overline{T\times
k^*(1,\ldots,1)}$.

The cone $K$ is strictly convex since $c_i\geq 0$ and  $\sum
c_i\chi_i'=0$ implies \ \ $c_i=0$\ \ for\ \ any\ \ $i$.\ \
Further,\ \ $\chi_i'$\  \ is \ \ an \ \ edge\ \ of\ \ $K$\ \
(i.e.,\ \ $\chi_i'\not\in
<\chi_1',\ldots,\chi'_{i-1},\chi'_{i+1},\ldots,\chi'_n>_{\Q_+}$)\
 if and only if $\chi_i$ is not contained in the convex hull of
the set $\{\chi_1,\ldots,\chi_{i-1},\chi_{i+1},\ldots,$
$\chi_n\}$. Finally, $\Q_+\chi_i'\ne \Q_+\chi_j'$ is equivalent to
$\chi_i\ne \chi_j$. \end{proof}

\begin{remark} If any hyperplane section of $X\subset \PProj$ is
reduced (i.e., it is a sum of prime divisors), then $X$ has (SP)
(see \cite[Lemma 2]{KrW}). If $X$ is the orbit of a highest vector
in an irreducible representation of a reductive group, then this
condition is equivalent to  (SP) (see \cite[Prop. 5]{KrW}). This
is not true for an orbit closure of a torus. Consider the action
of the torus $T=(k^*)^2 : \mathbb{P}(k^4)$,
$(t_1,t_2)(x_0:x_1:x_2:x_3)=(t_1t_2^2x_0, t_1t_2x_1, t_1^3x_2,
t_2^2x_3)$. The orbit closure of the point $(1:1:1:1)$ is the
hypersurface defined by the equation $x_0x_1^2=x_2x_3^2$. It has
(SP) and its intersection with $H_{x_0}$ is not reduced.
\end{remark}

\noindent {\bf Theorem 4'.} {\it The orbit closure
$X=\overline{T(1:\ldots:1)}$ satisfies (WSP) if and only if there
are no more then one $\chi_j$ in the interior of any face of the
convex hull ${conv}\{\chi_1,\ldots,\chi_n\}$ (in particular,
$\chi_i\ne\chi_j$ for $i\ne j$).}

\begin{proof} The proof is similar to the proof of Theorem 3'.
\end{proof}

\noindent {\bf Theorem 5'.} {\it The orbit closure
$X=\overline{T(1:\ldots:1)}\subset \PProj$ satisfies (SSP) if and
only if $X=\PProj$ (i.e., the weights $\chi_1,\ldots,\chi_n$ are
affinely independent).}

\begin{proof} The proof follows from Theorem 5.\end{proof}


\end{document}